\documentclass[12pt]{amsart}
\usepackage{amsmath, amssymb,amsthm, enumerate,graphicx}

\usepackage{fullpage}
\usepackage{setspace}

\newcommand{\mathd}{\mathrm{d}}
\newcommand{\mathi}{\mathrm{i}}
\newcommand{\nocomma}{}

\newcommand{\tmop}[1]{\ensuremath{\operatorname{#1}}}

\newcommand{\upl}{+}

\newenvironment{Proof}{\noindent\textbf{Proof\ }}{\hspace*{\fill}$\Box$\medskip}
\newenvironment{ProofB}{\noindent\textbf{Proof of Proposition B\ }}{\hspace*{\fill}$\Box$\medskip}
\newenvironment{ProofThm1}{\noindent\textbf{Proof of Theorem 1\ }}{\hspace*{\fill}$\Box$\medskip}
\newtheorem*{remark}{Remark}

\newtheorem{lemma}{Lemma}
\newtheorem*{lemma*}{Lemma}

\newtheorem{theorem}{Theorem}

\newtheorem*{propositionA}{Proposition A}
\newtheorem*{propositionB}{Proposition B}

\begin{document}
\title{Limitations to mollifying $\zeta(s)$.}
\author{Maksym Radziwi\l\l}
\address{Department of Mathematics \\ Stanford University \\
450 Serra Mall, Bldg. 380\\ Stanford, CA 94305-2125}
\email{maksym@stanford.edu}
\thanks{The author is partially supported by a NSERC PGS-D award.}
\begin{abstract}
 We establish limitations to how well one can mollify $\zeta(s)$ on the critical
line with mollifiers of arbitrary length. Our result gives a non-trivial lower bound for the contribution of the off-diagonal terms to mollified moments of $\zeta$. On the Riemann Hypothesis, we establish a connection between the mollified moment and Montgomery's Pair Correlation Function. 
\end{abstract}

\maketitle


\section{Introduction}\label{sec:intro}

The zero-distribution of an meromorphic function
and the distribution of its size are closely related problems
as can be seen from Jensen's inequality in complex analysis.
For this reason, when studying the zeros of 
the Riemann $\zeta$-function it is advantegeous to
reduce the size of $\zeta(s)$ and
to count instead the zeros of $\zeta(s)M(s)$ 
with $M(s)$ a \textit{mollifier}: an
entire function $M(s)$ pretending to behave
as $1 / \zeta(s)$ \cite{bl}. 
A natural choice for $M(s)$ is
$$
M(s) = \sum_{n \geq 1} \frac{\mu(n)W(n)}{n^s}
$$
with $W$ a smooth function ensuring the absolute convergence of the sum.

Away from the neighborhood of a zero of $\zeta(s)$, mollifiers are good 
pointwise approximations
to $1 / \zeta(s)$ (see \cite{Farmer}, Lemma 1)
Since there are at most a few
zeroes in the strip $\sigma > \tfrac 12 + \varepsilon$, a mollifier is
on average an excellent pointwise approximation to 
$1 / \zeta(s)$ to the right of the critical line.
On the critical line a mollifier is no longer a good pointwise approximation
to $1 / \zeta(s)$
because a positive proportion of the zeros lies on the critical line
\cite{SelZeroDensity}. 
For this reason on the half-line we consider
$$
\mathcal{I} = \mathcal{I}(M) := \frac{1}{T} \int_{T}^{2T} \big | 1 - \zeta(\tfrac 12 + \mathi t)M(\tfrac 12
+ \mathi t) \big |^2 \mathd t.
$$
The integral $\mathcal{I} $ is related to the horizontal
distribution of the zeros of $\zeta(s)$, for example via the inequality
$
\sum_{T \leq \gamma \leq 2T} |\beta - \tfrac 12| \ll T \log ( 1 + \mathcal{I}(M) ).
$ valid for any choice of Dirichlet polynomial $M$. 
Understanding $\mathcal{I}$, and in particular how small $\mathcal{I}(M)$ can be for various choices of $M$,
forms the principal focus of this paper. 

The mollifier
$$
\mathcal{L}_{\theta}(s) := \sum_{n \leq T^{\theta}} \frac{\mu(n)}{n^s} \cdot \bigg
(1 - \frac{\log n}{\log T^{\theta}} \bigg )
$$
achieves $\mathcal{I}(\mathcal{L}_{\theta}) \sim 1 / \theta$ for $\theta < \tfrac 47$ by a deep
result of Conrey \cite{ConreyPositive} (see also \cite{Balu}). 
It is
conjectured by Farmer \cite{Farmer} that, with this choice of
mollifier,
$\mathcal{I}(\mathcal{L}_{\theta}) \sim 1 / \theta$ for all $\theta > 0$. 
As we later show, among all Dirichlet polynomials
\begin{equation} \label{DirR}
M_{\theta}(s) = \sum_{n \leq T^{\theta}} \frac{a(n)}{n^s} \text{ with } a(n) \ll n^{\varepsilon}
\text{ and } a(1) = 1
\end{equation}
with $\theta < \tfrac 12$ fixed, the mollifier $\mathcal{L}_{\theta}(s)$ minimizes $\mathcal{I}$. 
We would like to understand if $\mathcal{I}(M_{\theta})$ can be much smaller 
than $1 / \theta$ when $M_{\theta}(s)$ is a longer mollifier,
say with $\theta > 1$. 
We show that the answer is ``no''. In fact, unconditionally, there is an
absolute
constant $c > 0$ such that $\mathcal{I}(M_{\theta}) \geq c / \theta$ for
all $\theta > 0$ and all $M_{\theta}$ as in (\ref{DirR}).
\begin{theorem} \label{Thm1}
Let $\theta > 0$ be given.
There is an absolute constant $c > 0$ such that
for all $T$ large enough, and  
all $M_{\theta}$ as in (\ref{DirR})
$$
\mathcal{I}(M_{\theta}) := \frac{1}{T}
\int_{T}^{2T} \big | 1 - \zeta(\tfrac 12 + \mathi t)M_{\theta}(\tfrac 12 + \mathi t)
\big |^2 \mathd t \geq \frac{c}{\theta}.
$$
\end{theorem}
The constant $c$ in Theorem 1 
depends on the proportion of the zeros of $\zeta(s)$
lying on the critical line.
The constant $c$ cannot be greater than one, since $c > 1$ would
contradict Farmer's conjecture in \cite{Farmer}. 
%
For $\theta < \tfrac 12$ we show that $c = 1$, using an asymptotic
formula for $\mathcal{I}$, due to Balasubramanian, Conrey and Heath-Brown.
Proposition B below is due to Prof. Soundararajan.
\begin{propositionB} [Soundararajan]
Let $M_{\theta}$ be as in (\ref{DirR}). If $\theta < \tfrac 12$, then, as $T \rightarrow \infty$, 
$$
\mathcal{I}(M_{\theta}) \sim 
\sum_{m,n \leq T^{\theta}} \frac{a(m)\overline{a(n)}}{[m,n]} \cdot \bigg
( \log \frac{T (m,n)^2}{2 \pi m n} + 2\log 2 + 2\gamma - 1 \bigg ) -1 
\geq \frac{1}{\theta} + o(1)
$$
\end{propositionB}
Similar quadratic forms have been considered by Selberg \cite{SelZeroDensity}
and Iwaniec-Sarnak \cite{IwaniecSarnak}. To the best of the authors knowledge this 
is the first time that the proof of such a lower bound appears in the
litterature. 

Proposition B suggests that most likely $c = 1$ for all $\theta > 0$. 
Assuming the
Riemann Hypothesis and the Pair Correlation conjecture we show that $c \geq 1 
- \varepsilon$ for all $\theta > \theta_0(\varepsilon)$ large enough.
This is interesting because one naively expects the problem to become more
difficult for large $\theta$.
\begin{theorem}
Let $\theta > 0$ be given. 
Assume the Riemann Hypothesis and the Pair Correlation Conjecture.
Let $M_{\theta}$ be as in (\ref{DirR}) and assume in addition that $a(p^k) \ll 1$.
Then, as $T \rightarrow \infty$,  
$$
\mathcal{I}(M_{\theta}) := \frac{1}{T} \int_{T}^{2T} \big | 1 - \zeta(\tfrac 12 + \mathi t)M_{\theta}(\tfrac 12 
+ \mathi t) \big |^2 \mathd t \geq \frac{1}{0.5 + \theta} \cdot 
( 1 + o_{\theta}(1)).
$$
\end{theorem}
\begin{remark}
The condition $a(p^k) \ll 1$ can be dispensed with.
\end{remark}
The size of $\mathcal{I}(M_{\theta})$
depends on the distribution of the zeros of $\zeta(s)$ in small interval of length
$2\pi/(1 + \theta)\log T$, around zeros of $\zeta(s)$. When $\theta$
is large, the Pair Correlation Conjecture allows to control the number of
zeros in such thin intervals, thus giving increasingly better
lower bounds for $\mathcal{I}(M_{\theta})$. 

On the Dirichlet polynomial side, an average of length $T$ such as 
in Theorem 2 detects the first $T$ coefficients of a Dirichlet series 
$F(s) = \sum a(n) n^{-s}$. 
This leads to a ``trivial'' lower bound (see \cite{RamBalu}),
$$
\int_{T}^{2T} |F(\tfrac 12 + \mathi t)|^2 \mathd t 
\gg 
T \sum_{n \leq T} \frac{|a(n)|^2}{n}.
$$
Let $F(s) = 1 - \zeta(s)M(s)$ with $M(s) = \sum_{n \leq T}
\mu(n) n^{-s}$. Then the first $T$ coefficients of $1 - \zeta(s)M(s)$ are zero
making the above lower bound vacuous. As another example let's consider
the Dirichlet series
$F(s) = 1 - \zeta(s)M(s)$ with $M(s) = \mathcal{L}_{\theta}(s)$. 
The trivial lower bound leads to
$c T / (1 + \theta)^2$ while Theorem 1 gives $c T / \theta$.

Theorems 1 and 2 beat the trivial lower bound by exploiting 
the relationship between $1 - \zeta(s)M(s)$
and the zeros of $\zeta(s)$. This is made explicit
in Proposition A below.
\begin{propositionA}
  Let $\varepsilon > 0$ and $\theta > 0$ be given. 
  Then for $T$ large,
  and for $S$ any $\delta := 2\pi A / \log T$ well-spaced
  subset of zeros of $\zeta(s)$
  with ordinates in $[T;2T]$, we have 
  for all $M_{\theta}$ as in (\ref{DirR})
%
  $$
  \frac{1}{T} \int_{T}^{2T} | 1 - \zeta(\tfrac 12 + \mathi t)M_{\theta}(\tfrac 12 
  + \mathi t)|^2 \mathd t \geq \frac{1 + O(\varepsilon)}{1 + \theta + \tfrac{1}{A}}
  \cdot \frac{\tmop{Card}(S)}{\tfrac{T}{2\pi} \log T}  + O(T^{\varepsilon}).
  $$
\end{propositionA}
The main idea in the proof of Proposition A is to connect, using
Sobolev's inequality, the
value of $1 - \zeta(s)M(s)$ at a zero with a continuous average of
$1 - \zeta(s)M(s)$ around that zero. 
Using this idea we can also give an elementary proof of a result of
Baez-Duarte, Balazard, Landreau and Saias \cite{Baez2}: For a mollifier
$M(s)$ of length $N$,
\begin{equation}
\label{Baezz}
\int_{\mathbb{R}} \bigg | \frac{1 - \zeta(\tfrac 12 + it)M(\tfrac 12 + it)}
{\tfrac 12 + it} \bigg |^{2} dt \geq \frac{C}{\log N}
\end{equation}
Their proof depends on functional analysis: by Plancherel (\ref{Baezz}) 
is related to the $L^2$ behavior of the function $\rho(x) = \{ 1 / x \}$.
Re-proving (\ref{Baezz}) was the starting point for this paper.

On the Riemann Hypothesis we obtain
an analogue of Proposition A 
involving Montgomery's Pair Correlation
function,
$$
F(\alpha, T) := \frac{2\pi}{T \log T} \sum_{T \leq \gamma, \gamma' 
\leq 2T} T^{\mathi \alpha  ( \gamma - \gamma')} \cdot w(\gamma - \gamma')
\text{ where } w(x) = \frac{4}{4 + x^2}.
$$
The function $F(\alpha, T)$
describes the vertical distribution of the zeros of $\zeta(s)$. 
Following Montgomery it is well known that $F(\alpha, T) = \alpha + o(1)$
for $\varepsilon \leq \alpha \leq 1$ and $F(\alpha, T) \geq o(1)$ for
all $\alpha$. The Pair
Correlation Conjecture is equivalent to 
$F(\alpha, T) = 1 + o(1)$ in $1 \leq \alpha 
\leq M$ for every fixed $M > 1$. Theorem 2 follows from
Theorem 3 below.
\begin{theorem}
Let $\theta > 0$ be given. Assume the Riemann Hypothesis. Let $M_{\theta}$ be as in (\ref{DirR}) and
assume in addition that $a(p^k) \ll 1$.
Then, for $T$ large,
$$
\mathcal{I}(M_{\theta}) = \frac{1}{T} \int_{T}^{2T} |1 - \zeta(\tfrac 12 + \mathi t)M_{\theta}(\tfrac 12 + \mathi t)|^2
\mathd t \geq \bigg ( \frac{1}{2} + \int_{1}^{1 + \theta
+ \varepsilon} F(\alpha, T) \mathd \alpha \bigg )^{-1}.
$$
\end{theorem}
\begin{remark}
As in Theorem 2 the requirement $a(p^k) \ll 1$ can be dispensed with
\end{remark}
In Theorem 3, choosing $M_{\theta}(s) = \mathcal{L}_{\theta}(s)$ for
$\theta < \tfrac 47$ and applying Conrey's result \cite{ConreyPositive}
we have $\mathcal{I}(\mathcal{L}_{\theta}) \sim \frac{1}{\theta}$ for $\tfrac 12 < \theta < \tfrac 47$ and thus, for $\tfrac 12 < \theta < \tfrac 47$, 
$$
\int_{1}^{1 + \theta} F(\alpha, T) \mathd \alpha > \theta - \tfrac 12 + o(1).
$$
as $T \rightarrow \infty$.
In a subsequent paper, we will improve this result assuming
the Generalized Riemann Hypothesis. Further we will investigate limitations
to mollifying $\zeta(s)$ in the context of Levinson's method. 

Theorems 1 and 2 have analogues
for double-mollifiers
$
M(s) =
\sum a(m,n) m^{-s} n^{-1 + s}
$.
In Theorem 1, for $\theta$ bounded away from zero, say $\theta > \tfrac{1}{100}$,  we can take
$
M(s) := \int \lambda^{-s} \mathd \mu(\lambda)
$
with
$\mu(\cdot)$ a finite measure, supported
in $[1;T^{\theta}]$ and such that
$\int_{1 \leq x \leq t} \mathd |\mu(x)| \ll t^{A}$ for some $A > 1$. 
In particular, for $\theta$ bounded away from zero,  the assumption $a(n) \ll n^\varepsilon$ 
in Theorem 1 can be relaxed to
$a(n) \ll n^{A}$ for some fixed $A > 0$.
\\
\\
\textbf{Acknowledgments.} I would like to thank my supervisor Kannan Soundararajan for his advice and encouragements, Bob Hough for many conversations 
concerning the subject of the paper and 
Sandro Bettin for a careful reading of this paper. 

\section{Key ideas}
Sobolev's inequality
\[ \left| f \left( x \right) \right| \leqslant \frac{1}{b - a} \int_a^b \left|
   f \left( u \right) \right| \mathd u + \int_a^b \left| f' \left( x \right)
   \right| \mathd x, \]
bounds a function $f$ at a particular point $a \leqslant x \leqslant b$, by an
average of $f$ and $f'$. For a Dirichlet polynomial $A \left( \cdot \right)$
we prove a Sobolev inequality without an average over $A'$.

\begin{lemma}
  \label{lemma1}Let $A$ be a Dirichlet polynomial supported on integers $n$
  with $M \leqslant n \leqslant N$. If $f$ is a smooth function such that $f
  \left( x \right) = 1$ for $\log M \leqslant 2 \pi x \leqslant \log N$, then
  for all real $u$,
  \[ A \left( \mathi u \right) = \int_{- \infty}^{\infty} A \left( \mathi t
     \right) \hat{f} \left( t - u \right) \mathd t . \]
\end{lemma}

\begin{Proof}
  Expanding $A \left( s \right) = \sum_{M \leqslant n \leqslant N} a \left( n
  \right) n^{- s}$ and using Fourier inversion,
  \begin{eqnarray*}
    \int_{- \infty}^{\infty} A \left( \mathi t \right) \hat{f} \left( t - u
    \right) \mathd t & = & \sum_{M \leqslant n \leqslant N} a \left( n \right)
    \int_{- \infty}^{\infty} n^{- \mathi t} \cdot \hat{f} \left( t - u \right)
    \mathd t\\
    & = & \sum_{M \leqslant n \leqslant N} a \left( n \right) n^{- \mathi u}
    \cdot f \left( \frac{\log n}{2 \pi} \right) .
  \end{eqnarray*}
  By assumptions, $f \left( \log n / (2 \pi) \right) = 1$ for $M \leqslant n
  \leqslant N$, and so the right-hand side is equal to $A \left( \mathi u
  \right)$.
\end{Proof}

In the above lemma we can take $\zeta \left( s \right)$ or $1 - \zeta \left( s
\right) A \left( s \right)$ instead of $A \left( s \right)$ because $\zeta
\left( s \right)$ is approximated very well by a Dirichlet polynomial.

\begin{lemma}
  \label{lemma2}There is a smooth function $w \left( x \right)$ with $0
  \leqslant w \left( x \right) \leqslant 1$, $w \left( 0 \right) = 1$, such
  that for $T \leqslant t \leqslant 2 T$ , $T_1 = T^{1 \upl \varepsilon}$, and
  any fixed $v > 0$,
  \[ \zeta \left( s \right) = \sum_{n \leqslant T_1} n^{- s} \cdot w \left(
     \frac{n}{T_1} \right) + O_v \left( T^{- v} \right) . \]
\end{lemma}

\begin{Proof}
  This is Proposition 1 in Bombieri-Friedlander {\cite{Bombieri}}.
\end{Proof}

If $M$ is a long mollifier and $s$ is away from a zero of $\zeta (s)$ (on a
scale of $2 \pi / \log |s|$) then $1 - \zeta \left( s \right) M \left( s
\right) \approx 0$ . On the other hand, if on the same scale $s$ is close to a
zero of $\zeta (s)$, then $\zeta (s) M (s) \approx 0$ and therefore $1 - \zeta
(s)M (s) \approx 1$. Given a smooth $\hat{f} (x)$ concentrated in $|x| \ll 2
\pi / \log T$, the function
\begin{equation}
  \label{emulator} \sum_{\rho} \hat{f} \left( t - \gamma \right) \text{,} T
  \leq t \leq 2 T
\end{equation}
exhibits a similar behavior to that of $1 - \zeta (s) M (s)$. However,
understanding the mean-square of (\ref{emulator}) is much simpler.

\begin{lemma}
  \label{lemma3}Let $S$ be a finite set and $f$ be a smooth function. If $K$
  is a smooth function with $K \geqslant f^2$, then,
  \[ \int_{- \infty}^{\infty} \left| \sum_{\gamma \in S} \hat{f} \left( t -
     \gamma \right) \right|^2 \mathd t \leqslant \sum_{\gamma, \gamma'}
     \hat{K} \left( \gamma - \gamma' \right) . \]
\end{lemma}

\begin{Proof}
  Notice that,
  \begin{eqnarray*}
    \sum_{\gamma \in S} \hat{f} \left( t - \gamma \right) & = & \sum_{\gamma
    \in S} \int_{- \infty}^{\infty} f \left( v \right) \cdot e^{2 \pi \mathi v
    \left( t - \gamma \right)} \mathd v\\
    & = & \int_{- \infty}^{\infty} e^{2 \pi \mathi vt} \cdot f \left( v
    \right) \sum_{\gamma \in S} e^{- 2 \pi \mathi \gamma v} \mathd v .
  \end{eqnarray*}
  Therefore by Plancherel,
  \begin{eqnarray*}
    \int_{- \infty}^{\infty} \left| \sum_{\gamma \in S} \hat{f} \left( t -
    \gamma \right) \right|^2 \mathd t & = & \int_{- \infty}^{\infty} \left|
    \sum_{\gamma \in S} e^{- 2 \pi \mathi \gamma v} \right|^2 \cdot \left| f
    \left( v \right) \right|^2 \mathd v\\
    & \leqslant & \int_{- \infty}^{\infty} \left| \sum_{\gamma \in S} e^{- 2
    \pi \mathi \gamma v} \right|^2 \cdot K \left( v \right) \mathd v .
  \end{eqnarray*}
  Expanding the square, we find
  \[ \sum_{\gamma, \gamma' \in S} \int_{- \infty}^{\infty} e^{2 \pi \mathi v
     \left( \gamma - \gamma' \right)} \cdot K \left( v \right) \mathd v =
     \sum_{\gamma, \gamma' \in S} \hat{K} \left( \gamma - \gamma' \right), \]
  as desired.
\end{Proof}

For a $\delta$-well-spaced set $S$ it is convenient to pick a $K$ such that
$\hat{K} \left( x \right) = 0$ when $\left| x \right| \geqslant \delta$. For
such a choice of $K$,
\[ \sum_{\gamma, \gamma' \in S} \hat{K} \left( \gamma - \gamma' \right) =
   \hat{K} \left( 0 \right) \cdot \tmop{Card} \left( S \right) . \]
We construct in the lemma below a set of functions with this property. These
are known as the Beurling-Selberg majorants.

\begin{lemma}
  \label{lemma4}\label{Beurling}Let $\delta > 0$.
  For any interval $I = [a,b]$,
  there exists an even entire
  function $K \left( w \right)$ such that,
  \begin{itemize}
    \item $K \left( u \right) \geqslant \chi_I \left( u \right)$
    
    \item $\hat{K} \left( 0 \right) = b - a + 1 / \delta$
    
    \item $\hat{K} \left( x \right) = 0$ for $\left| x \right| > \delta$.
  \end{itemize}
\end{lemma}

\begin{Proof}
  Beurling {\cite{Beurling2}} considered the function,
  \[ B \left( z \right) = \left( \frac{\sin \pi z}{\pi} \right)^2 \cdot \left(
     \frac{1}{z^2} + \sum_{n = 0}^{\infty} \frac{1}{\left( z - n \right)^2} -
     \sum_{n = 1}^{\infty} \frac{1}{\left( z + n \right)^2} \right) . \]
  The function $B \left( z \right)$ is entire, has the property that $B \left(
  x \right) \geqslant \tmop{sgn} \left( x \right)$, and
  \[ \int_{- \infty}^{\infty} B \left( x \right) - \tmop{sgn} \left( x \right)
     \mathd x = 1. \]
  From the definition of $B \left( z \right)$ it is easy to see that $B \left(
  z \right) = O \left( e^{2 \pi \left| \tmop{Imz} \right|} \right)$.
  Therefore, by Paley-Wiener $\hat{B} \left( x \right) = 0$ when $\left| x
  \right| \geqslant 1$. Given an interval $I = \left[ a, b \right]$ we define
  \[ K \left( z \right) = \frac{1}{2} \cdot B \left( \delta \left( z - a
     \right) \right) + \frac{1}{2} \cdot B \left( \delta \left( b - z \right)
     \right) . \]
  Then by a direct check using the properties of $B \left( z \right)$ we find
  that, $K \left( x \right) \geqslant \chi_I \left( x \right)$ for all real
  $x$, $\hat{K} \left( x \right) = 0$ for $\left| x \right| \geqslant \delta$,
  and $\hat{K} \left( 0 \right) = \int_{\mathbb{R}} K \left( x \right) \mathd
  x = b - a + 1 / \delta$, as desired.
\end{Proof}

\section{Proof of Proposition A.}

We denote by $t$ the imaginary part of $s$. Let $\eta > 0$. By Lemma
\ref{lemma2} there is a smooth function $w \left( x \right)$ with $0 \leqslant
w \left( x \right) \leqslant 1$, $w \left( 0 \right) = 1$, and such that for
$T \leq t \leq 2 T$,
\[ \zeta \left( s \right) = \sum_{n \leqslant T^{1 + \eta}} \frac{1}{n^s}
   \cdot w \left( \frac{n}{T^{1 + \eta}} \right) + O_v \left( T^{- v} \right)
   . \]
Multiplying by a Dirichlet polynomial $M \left( s \right) = \sum a \left( n
\right) n^{- s}$ of length $N = T^{\theta}$ and with coefficients bounded by
$N$ we obtain a Dirichlet polynomial $B (s)$ of length $T^{1 + \eta} \cdot N =
T^{1 + \eta + \theta}$ for which,
\begin{equation}
  \label{mastereq00} \zeta \left( s \right) M \left( s \right) = B \left( s
  \right) + O_v \left( T^{- v} \right) .
\end{equation}
Since $1 - B \left( s \right)$ is a Dirichlet polynomial of length $T^{1 +
\eta} \cdot N$, by Lemma \ref{lemma1}, for any smooth function $f$ with $f
\left( x \right) = 1$ in $1 \leqslant 2 \pi x \leqslant \log (T^{1 + \eta}
\cdot N)$,
\begin{equation}
  \label{mastereq01} 1 - B \left( \tfrac{1}{2} + \mathi u \right) = \int_{-
  \infty}^{\infty} \left( 1 - B \left( \tfrac{1}{2} + \mathi t \right) \right)
  \hat{f} \left( t - u \right) \mathd t .
\end{equation}
We choose a function $f$ supported on the interval $0 \leq 2 \pi x \leq \log (T^{1 + \eta} \cdot N) +
1$, equal to one in $1 \leq 2 \pi x \leq \log (T^{1 + \eta} \cdot N)$ and
bounded between $0$ and $1$, with $f^{(\ell)} (x) \ll_{\ell} 1$ for any given
$\ell > 0$. Here is an example of such a function $f$,
\begin{center}
  \includegraphics[scale=0.6]{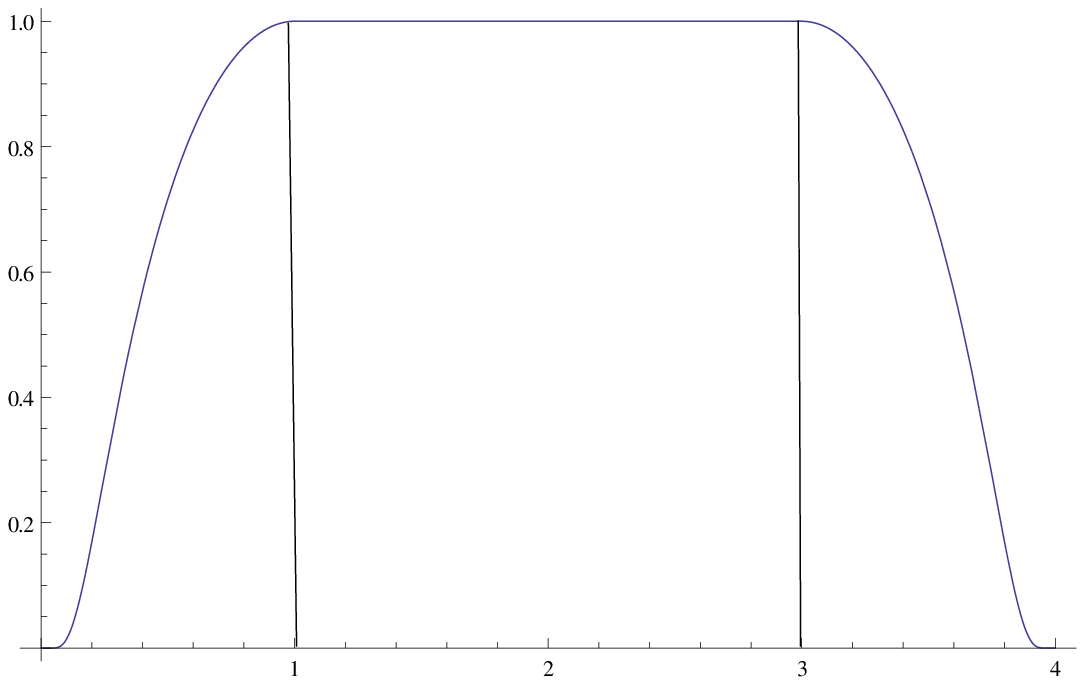}
\end{center}
For any fixed $v$, $\hat{f} \left( x \right) \ll \left( \log T \right) \cdot
\left( 1 + \left| x \right| \log T \right)^{- v}$. Therefore for $T + T^{\eta}
\leqslant u \leqslant 2 T - T^{\eta}$ and $t \notin [T ; 2 T]$, we have
$\hat{f} (t - u) \ll_v T^{- \eta v / 2} \cdot (1 + |x - u| \log T)^{- v / 2}
\ll_{\eta, v} T^{- v} \cdot (1 + |x - u| \log T)^{- v}$. 
Since $1 - B ( \tfrac{1}{2}
+ \mathi t)$ is trivially bounded by $N^2 \ll T^{2 \theta}$ we get
for $T + T^{\eta} \leq u \leq 2T - T^{\eta}$,
\begin{equation}
  \label{equation2} \int_{- \infty}^{\infty} (1 - B ( \tfrac{1}{2} + \mathi
  t)) \hat{f} (t - u) \mathd t = \int_T^{2 T} (1 - B ( \tfrac{1}{2} + \mathi
  t)) \hat{f} (t - u) \mathd t + O_{\eta, v} (T^{- v}) .
\end{equation}
Combining (\ref{mastereq01}) with (\ref{equation2}) and (\ref{mastereq00}) we
obtain
\[ 1 - \zeta \left( \tfrac{1}{2} + \mathi u \right) M \left( \tfrac{1}{2} +
   \mathi u \right) = \int_T^{2 T} \left( 1 - \zeta \left( \tfrac{1}{2} +
   \mathi t \right) M \left( \tfrac{1}{2} + \mathi t \right) \right) \hat{f}
   \left( t - u \right) \mathd t + O_{\eta, v} \left( T^{- v} \right) . \]
In the above equation take $u = \gamma$, with $\gamma$ the ordinate of a zero
of $\zeta \left( s \right)$ lying on the half-line and with $T + T^{\eta}
\leqslant \gamma \leqslant 2 T - T^{\eta}$. Summing over any set $S$ of such
zeros, we get
\[ \tmop{Card} \left( S \right) = \int_T^{2 T} \left( 1 - \zeta \left(
   \tfrac{1}{2} + \mathi t \right) M \left( \tfrac{1}{2} + \mathi t \right)
   \right) \sum_{\gamma \in S} \hat{f} \left( t - \gamma \right) \mathd t +
   O_{\eta, v} \left( T^{- v} \right) . \]
By Cauchy-Schwarz
\begin{multline*} \text{\tmop{Card}} \left( S \right) \leqslant \bigg 
( \int_T^{2 T} \bigg | 1
   - \zeta \left( \tfrac{1}{2} + \mathi t \right) M \left( \tfrac{1}{2} +
   \mathi t \right) \bigg |^2 \mathd t \bigg )^{1 / 2} \cdot \bigg (
   \int_{\mathbb{R}} \bigg | \sum_{\gamma \in S} \hat{f} \left ( t - \gamma
   \right) \bigg|^2 \mathd t \bigg )^{1 / 2}  \\ + O_{\eta, v} \left( T^{- v}
   \right) .
\end{multline*}
By Lemma \ref{lemma3}, for any $K$ such that $K \geqslant f^2$,
\[ \int_{\mathbb{R}} \left| \sum_{\gamma \in S} \hat{f} \left( t - \gamma
   \right) \right|^2 \mathd t \leqslant \sum_{\gamma, \gamma' \in S} \hat{K}
   \left( \gamma - \gamma' \right) . \]
Since $0 \leqslant f \leqslant 1$ and $f$ is supported in $I = \left[ 0 ;
\left( 1 / 2 \pi \right) \cdot \log \left( eT^{1 + \eta} \cdot N \right)
\right]$ the condition $K \geqslant f^2$ is satisfied whenever $K \geqslant
\chi_I$. Using Lemma \ref{lemma4}, we pick a function $K$ such that, $K
\geqslant \chi_I$, $\hat{K} \left( x \right) = 0$ for $\left| x \right|
\geqslant \delta := 2\pi A / \log T$, 
and $\hat{K} \left( 0 \right) = \left| I \right| + 1 /
\delta$. Since the set $S$ is $\delta$ well-spaced,
\begin{align*}
  \sum_{\gamma, \gamma' \in S} \hat{K} \left(\gamma - \gamma' \right) &
  =\hat{K}\left(0\right){\cdot} \text{Card} \left(S\right)=\left(|I|+1/{\delta}\right) \cdot {\text{Card}}(S)\\
  & = (1+O({\eta})) \cdot 
  \frac{\log T}{2\pi} \cdot
  (1 + \theta + \tfrac{1}{A}) \cdot \text{Card}(S).
\end{align*}
Combining the above three equations, we conclude
\[ T \cdot \frac{1 + O(\eta)}{1 + \theta + \tfrac{1}{A}} 
\cdot \frac{\text{Card(S)}}{\tfrac{T}{2\pi} \log T} \leqslant
   \int_T^{2 T} \left| 1 - \zeta \left( \tfrac{1}{2} + \mathi t \right) M
   \left( \tfrac{1}{2} + \mathi t \right) \right|^2 \mathd t + O_{\eta, v}
   (T^{- v}) . \]
At the price of an additional error term $O (T^{\eta} \cdot \log T)$ we can
add to $S$ an arbitrary set of zeros with ordinates $\gamma$ in the interval
$\left[ T ; T + T^{\eta} \right] \cup \left[ 2 T - T^{\eta} ; 2T^{\eta}
\right]$. Taking $\eta \rightarrow 0$ very slowly as $T \rightarrow \infty$ we
obtain the claim.

\section{Deduction of Theorem 1}

Theorem \ref{Thm1} follows from Proposition A and the existence of a well
spaced set of zeros, lying on the critical line, with ordinates in $[T;2T]$
and cardinality $\gg N
\left( T \right) \asymp T \log T$.

\begin{lemma}
  \label{zero}There is a set $S$ of zeros of $\zeta \left( s \right)$ with
  $\beta = \tfrac{1}{2}$ and $T \leqslant \gamma \leqslant 2 T$, such that
  \begin{itemize}
    \item The elements of $S$ are $2\pi A / \log T$ well-spaced, for some absolute
    constant $A > 0$.
    
    \item The set $S$ has $\gg T \log T$ elements.
  \end{itemize}
\end{lemma}

\begin{Proof}
  Selberg's proof (\cite{Titchmarsh}, 10.22, p. 279)
  shows that there is an $h = 2\pi A / \log T$,
  with $A > 0$ constant, for which the set
  \[ E = \left\{ T \leqslant t \leqslant 2 T : \gamma \in \left( t ; t + h
     \right) \text{ for some } \rho = \tfrac{1}{2} + \mathi \gamma
     \right\}, \]
  has $\tmop{meas} \left\{ E \right\} \geqslant c \cdot T$ with $c > 0$
  constant. Hence at least $c \cdot T / h$ intervals $(T + nh ; T + (n + 1)
  h)$ contain a $t$ such that there is a zero with $\beta = \tfrac{1}{2}$ and
  $\gamma \in (t ; t + h)$. It follows that at least $c \cdot T / (2 h)$
  intervals $(T + (n - 1) h ; T + (n + 1) h)$ contain the ordinate of a zero
  lying on the half-line. Taking every third such intervals produces a
  sequence of $c \cdot T / 6 h$ intervals of length $2 h$, and spaced by at
  least $h$, each containing the ordinate of a zero on the half-line. Thus we
  obtain a $h$ well-spaced set $S$ of at least $\geqslant c \cdot T / 6 h$
  zeros of $\zeta (s)$ lying on the half-line, with ordinates in $T \leqslant
  \gamma \leqslant 2 T$.
\end{Proof}

\begin{ProofThm1}
  By Proposition A, given $\varepsilon > 0$,
  for any $2\pi A / \log T$-well spaced
  set of zeros $S$ of $\zeta \left( s \right)$ lying on the critical line and with
  ordinates in $[T;2T]$,
  \begin{equation} \label{liminf}
    \frac{1}{T}\int_T^{2 T} \left| 1 - \zeta \left( \tfrac{1}{2} + \mathi t \right) M_{\theta}
     \left( \tfrac 12 + \mathi t \right) \right|^2 \mathd t \geqslant 
     \frac{\text{Card}(S)}{\tfrac{T}{2\pi}\log T} \cdot 
     \frac{1 + O(\varepsilon)}{1 + \theta + 1/A} 
   \end{equation}
  We pick $S$ as in Lemma \ref{zero}. Then, for $\theta > \tfrac{1}{2}$ the
  above lower bound is,
  $$
  \geq c_1 \frac{(1 + O(\varepsilon))}{1 + \theta} \geq c_2 \frac{1 + O(\varepsilon)}{\theta}
  $$
  with $c_1, c_2 > 0$ absolute constants. 
  Since $\varepsilon > 0$ is arbitrary, it follows that the $\text{liminf}$ of the
  left-hand side of (\ref{liminf}) is at least $c / \theta$, as desired. 
  On the other hand when $\theta < \tfrac 12$,
  Theorem \ref{Thm1} follows from Proposition B. 
\end{ProofThm1}

\section{Preliminaries for Theorem 2 and 3}

The proof of Theorem 2
follows the lines of proof of Proposition A. There are
two main differences. The first is that for $n \ll T^{1 - \varepsilon}$ we
exploit cancellations in the sum,
\[ \label{Gonek} \sum_{T \leqslant \gamma \leqslant 2 T} n^{- \mathi \gamma} .
\]
This is possible because we assume the Riemann Hypothesis.

\begin{lemma}
  \label{lemma6}Assume the Riemann Hypothesis. Uniformly in integer $n
  \geqslant 2$,
  \[ \sum_{T \leqslant \gamma \leqslant 2 T} n^{-1/2 - \mathi \gamma} = -
     \frac{T}{2 \pi} \cdot \frac{\Lambda \left( n \right)}{n} + O \left(
     \left( \log T \right)^2 \cdot n \right) . \]
\end{lemma}

\begin{Proof}
  See Gonek's paper {\cite{GonekLemma}}.
\end{Proof}

\begin{lemma} \label{lemma7}
  Let $A \left( s \right) = \sum a \left( n \right) \cdot n^{- s}$ be a
  Dirichlet polynomial of length $N$. Let $f$ be a smooth test function. Then,
  for real $u$,
  \[ \int_{- \infty}^{\infty} A \left( \mathi u \right) \hat{f} \left( t -
     \gamma \right) \mathd t = \sum_{n \leqslant M} \frac{a \left( n
     \right)}{n^{\mathi u}} \cdot f \left( \frac{\log n}{2 \pi} \right). \]
\end{lemma}

\begin{Proof}
  Expanding $A \left( s \right) = \sum a \left( n \right) \cdot n^{- s}$ and
  using Fourier inversion,
  \begin{eqnarray*}
    \int_{- \infty}^{\infty} A \left( \mathi t \right) \hat{f} \left( t - u
    \right) \mathd t & = & \sum_{n \leqslant N} a \left( n \right) \int_{-
    \infty}^{\infty} n^{- \mathi t} \cdot \hat{f} \left( t - u \right) \mathd
    t\\
    & = & \sum_{n \leqslant N} a \left( n \right) n^{- \mathi u} \cdot f
    \left( \frac{\log n}{2 \pi} \right).
  \end{eqnarray*}
  as claimed. 
\end{Proof}

The second difference with the proof of Proposition A, is that on the Riemann Hypothesis 
we can estimate asymptotically sums of the form
\[ \sum_{T \leqslant \gamma, \gamma' \leqslant 2 T} \hat{K} \left( \gamma -
   \gamma' \right) . \]
In application $\hat{K} \left( x \right)$ will be concentrated in $\left| x
\right| \ll 1 / \log T$, so that by the uncertainty principle, $K \left( x
\right)$ will be spread out on intervals of length $\asymp \log T$ (or
longer). If the Pair Correlation conjecture is not assumed then
the lemma below is true with Montgomery's Pair Correlation $F(\alpha, T)$
instead of its limit $F(\alpha)$. 

\begin{lemma}
  \label{PairCorel}Assume the Riemann Hypothesis.
  Let $h \geqslant 0$ denote a smooth, non-zero, and compactly
  supported function. Let $K \left( x \right) = h \left( 2 \pi x / \log T
  \right)$. Then, as $T \rightarrow \infty$,
  \[ \sum_{T + T^{\varepsilon} \leqslant \gamma, \gamma' \leqslant 2 T -
     T^{\varepsilon}} \hat{K} \left( \gamma - \gamma' \right) = T \cdot
     \left( \frac{\log T}{2 \pi} \right)^2 \int_{- \infty}^{\infty} h \left(
     \alpha \right) \cdot F \left( \alpha, T \right) \mathd x
     + O(T^{1 - \varepsilon})\]
  with $F(\alpha, T)$ Montgomery's Pair Correlation function.
\end{lemma}

\begin{Proof}
  Since $K \left( x \right) = K \left( \log T / 2 \pi \cdot x \right)$ the
  Fourier transform of $K$ is given by,
  \[ \hat{K} \left( x \right) = \frac{\log T}{2 \pi} \cdot \hat{h} \left(
     \frac{\log T}{2 \pi} \cdot x \right) . \]
  By definition
  \begin{align*}
    \sum_{T{\leqslant}{\gamma}, {\gamma}'{\leqslant}2T}\hat{h}\left({\frac{\log
    T}{2{\pi}}}{\cdot}({\gamma}-{\gamma}')\right) w(\gamma - \gamma')
    & ={\frac{T \cdot \log
    T}{2{\pi}}}\int_{-{\infty}}^{{\infty}}h\left({\alpha}\right)F\left({\alpha}, T \right){\mathd}{\alpha}.
  \end{align*}
with the weight $w(x) = 4 / (4 + x^2)$. 
  Multiplying by $\log T / 2\pi$, we obtain,
  \begin{equation} 
    \sum_{T \leqslant \gamma, \gamma' \leqslant 2 T} \hat{K} \left( \gamma -
    \gamma' \right) w(\gamma - \gamma') \sim T \cdot \left( \frac{\log T}{2 \pi} \right)^2 \int_{-
    \infty}^{\infty} h \left( \alpha \right) F \left( \alpha \right) \mathd
    \alpha . \label{paircorel222}
  \end{equation}
  One removes the weight $w(\gamma - \gamma')$ by a standard argument
  which we omit.
  Since $h$ is smooth, and compactly supported we have $\hat{K} \left( x
  \right) \ll_v \left( \log T \right) \cdot \left( 1 + \log T \cdot \left| x
  \right| \right)^{- v}$ for any fixed $v$. Thus, for any $\gamma,$
  \[ \sum_{T \leqslant \gamma \leqslant 2 T} \hat{K} \left( \gamma - \gamma'
     \right) \ll \left( \log T \right)^2 \]
  Since there are at most $\ll T^{\varepsilon} \cdot \log T$ ordinates of
  zeros in $\left[ T ; T + T^{\varepsilon} \right] \cup \left[ 2 T -
  T^{\varepsilon} ; 2 T \right]$, we can restrict the summation in
  (\ref{paircorel222}) to $T + T^{\varepsilon} \leqslant \gamma, \gamma'
  \leqslant 2 T - T^{\varepsilon}$ at the price of a negligible error term
  $\ll T^{\varepsilon} \cdot \left( \log T \right)^3$.
\end{Proof}

\section{Proof of Theorem 2 and 3}

We denote by $t$ the imaginary part of $s$.
Let $M$ be a Dirichlet polynomial of length $N = T^{\theta}$. 
Fix a small $
\tfrac {1}{10} > \eta > 0$. 
Proceeding as in the proof of Proposition A, there is a Dirichlet
polynomial
$B(s)$ of length $T^{1 + \eta} N$ such that for $T \leq t \leq
2T$ and for any fixed $v > 0$, 
\begin{equation} \label{mastereq}
  \zeta \left( s \right) M \left( s \right) = B \left( s \right) + O_v \left(
  T^{- v} \right).
\end{equation}
Since $a(1) = 1$, $a(p^k) \ll 1$ and $a(n) \ll n^{\varepsilon}$, the coefficients $b(n)$ of $B(s)$ satisfy, 
\[ b \left( 1 \right) = 1 + O (T^{- 1 - \eta}) \nocomma, \text{ } b
   (p^k) \ll 1, \text{ and } b \left( n \right) \ll n^{\varepsilon}. \]
Let $h(x) = h_0(2\pi x / \log T)$ with $h_0 \leq 1$ a smooth function
supported on $[\eta ; 1 + \theta + 2\eta]$ 
and equal to one on $[2\eta; 1 + \theta + \eta]$. 
These requirements on $h$ force that $\hat{h}(x) \ll_{\ell} \log T 
\cdot ( 1 + \log T |x|)^{-\ell}$ for every fixed $\ell > 0$. 

\begin{lemma*} 
We have
\begin{equation} \label{importanteq}
  \int_{- \infty}^{\infty} \left( 1 - B \left( \tfrac 12  + \mathi t \right)
  \right) \sum_{\gamma \in S} \hat{h} \left( t - \gamma \right) \mathd t =
  \left( 1 + O \left( \eta \right) \right) N \left( T \right).
\end{equation}

\end{lemma*}
\begin{Proof}
Write $h = f - g$ with $f(x)=f_0(2\pi x / \log T),g(x) = g_0(2\pi x / \log T)$
two smooth compactly supported functions such that $f_0(x ) = 1$ on $[0; 1 + \theta + \eta]$,
$g_0(x) = 1$ on $[0;\eta]$ and $g_0(x)$ is supported
on $[-A; 2\eta]$ for some $A > 0$. 
 By Lemma
\ref{lemma7} applied to $1 - B(s)$, 
\[ \int_{- \infty}^{\infty} \left( 1 - B \left( \tfrac 12 + \mathi t \right)
   \right) \hat{g} \left( t - u \right) \mathd t = 1 - b \left( 1 \right) +
   \sum_{2 \leqslant n \leqslant T^{2\eta}} \frac{b \left( n
   \right)}{n^{1/2 + \mathi u}} \cdot g \big( \frac{\log n}{2 \pi} \big).
\]
Set $u = \gamma$, and sum over the set $S$ of all zeros with ordinates $T +
T^{\eta} \leqslant \gamma \leqslant 2 T - T^{\eta}$. Using Gonek's
Lemma \ref{lemma6} and $1 - b(1) \ll T^{-1 - \eta}$, $b(p^k) \ll
1$, $g \ll 1$, we
get
\begin{align} \nonumber
\int_{- \infty}^{\infty} \left( 1 - B \left( \tfrac 12  + \mathi t \right)
   \right) \sum_{\gamma \in S} \hat{g} \left( t - \gamma \right) \mathd t
   & = -
   \frac{T}{2 \pi} \sum_{n \leqslant T^{2\eta}} \frac{b \left( n \right)
   \Lambda \left( n \right)}{n} \cdot g \big( \frac{\log n}{2 \pi} \big) + O
   \left( T^{3\eta} \right) \\ \label{final1}
   & \ll T \sum_{n \leq T^{2\eta}} \frac{\Lambda(n)}{n} \ll
   \eta T \log T \ll \eta N(T).
\end{align}
Since $1 - B(s)$ is of length $T^{1 + \eta}
N$, and $f(x) = 1$ on $1 \leq 2\pi x \leq \log(T^{1 +
  \eta} N)$; we get by Lemma \ref{lemma1},
\begin{equation*}
  \int_{- \infty}^{\infty} \left( 1 - B \left( \tfrac 12  + \mathi t \right)
  \right) \hat{f} \left( t - u \right) \mathd t = 1 - B(\tfrac 12 + \mathi
  u) + O_v \left( T^{-v} \right).
\end{equation*}
Set $u = \gamma$ and note that by equation (\ref{mastereq}), $B(\tfrac 12 + \mathi \gamma) = O_v(T^{-v})$. Summing over all $T + T^{\eta} \leq \gamma
\leq 2T - T^{\eta}$ we obtain
\begin{equation} \label{final2}
  \int_{- \infty}^{\infty} \left( 1 - B \left( \tfrac 12  + \mathi t \right)
  \right) \sum_{\gamma \in S} \hat{f} \left( t - \gamma \right) \mathd t =
  N(T) + O(T^{\eta}).
\end{equation}
Subtracting (\ref{final1}) from (\ref{final2}), 
and recalling that $\hat{h} = \hat{f} - \hat{g}$ (because $h = f - g$), 
we obtain the claim.  \end{Proof}

Since $\hat{h}(x) \ll_{v} \log T \cdot ( 1 + \log T |x| )^{-v}$, we proceed
exactly as in 
the proof of Proposition A; we truncate the integral
in (\ref{importanteq}) at $T$ and $2 T$, and using (\ref{mastereq}) 
replace $1 - B \left( s \right)$
by $1 - \zeta \left( s \right) M \left( s \right)$. Thus we obtain from 
(\ref{importanteq}) that
\[ \int_T^{2 T} \left( 1 - \zeta \left(\tfrac 12  + \mathi t \right) M \left( 
\tfrac 12 + \mathi t \right) \right) \sum_{\gamma \in S} \hat{h} \left( t - \gamma
   \right) \mathd t = \left( 1 + O \left( \eta \right) \right) N \left(
   T \right). \]
Applying Cauchy-Schwarz leads to
\[ ( 1 + O(\eta)) N(T) \leqslant \left( \int_T^{2 T} |1 - \zeta \left(\tfrac 12  + \mathi t \right) M   \left(\tfrac 12  + \mathi t \right) |^2 \mathd t \right)^{1 / 2} \cdot \left(
   \int_{\mathbb{R}} \bigg| \sum_{\gamma \in S} \hat{h} \left( t - \gamma
   \right) \bigg |^2 \mathd t \right)^{1 / 2}. \]
By Lemma \ref{lemma3}, for any smooth $K$ with $K \geqslant h^2$,
\begin{eqnarray*}
  \int_{\mathbb{R}} \bigg| \sum_{\gamma \in S} \hat{h} \left( t - \gamma
  \right) \bigg|^2 \mathd t & \leqslant & \sum_{\gamma, \gamma' \in S}
  \hat{K} \left( \gamma - \gamma' \right).
\end{eqnarray*}
Take $K = h^2$, and recall that $h(x) = h_0(2 \pi x / \log T) \leq 1$ 
with $h_0$ supported on $[\eta ; 1 + \theta + 2\eta]$. 
Thus $\hat{K}(x) = \log T / 2\pi \cdot \hat{h_0^2}(x \log T / 2\pi)$.  
Applying Lemma \ref{PairCorel} and bounding $h_0$ by 
$1$ on its interval of support we obtain
\begin{align*}
  \sum_{\gamma, \gamma' \in S} \hat{K} \left( \gamma - \gamma' \right) 
  & = \frac{\log T}{2\pi} \sum_{\gamma, \gamma' \in S} \hat{h_0^2} 
  \left ( \frac{\log T}{2\pi} (\gamma - \gamma') \right ) \\
  & \sim
  \frac{T \left (\log T \right )^2}{(2 \pi)^2} 
  \int_{- \infty}^{\infty} h_0^2 \left( \alpha \right) F
  \left( \alpha, T \right) \mathd \alpha\\
  & \leqslant ( 1 + o(1))
  \frac{T \left( \log T \right)^2}{(2 \pi)^2} \int_{\eta}^{1 \upl \theta +
  2\eta} F \left( \alpha, T \right) \mathd \alpha .
\end{align*}
Finally recall that $F(\alpha, T) = |\alpha| + o(1)$ uniformly for
$\eta < |\alpha| < 1 - \eta$. Combining the above three
inequalities and letting $\eta \rightarrow 0$ 
very slowly as $T \rightarrow \infty$ establishes Theorem 3.

To prove Theorem 2 note that on the Pair Correlation
conjecture (PCC) $F(\alpha, T) = 1 + o(1)$ for $1 \leq |\alpha| \leq M$,
and any fixed $M > 1$. Therefore on PCC,
$$
\int_{\eta}^{1 + \theta + 2\eta} F(\alpha, T) d\alpha = 
0.5 + \theta + O(\eta)
$$
Combining the above four equations and letting $\eta \rightarrow 
0$ we obtain Theorem 2. Alternatively, Theorem 2 is
an immediate consequence of Theorem 3 as explained in the
introduction.
%
%

\section{Proof of Proposition B.}

Our goal is to determine the minimum of the quadratic form,
\begin{equation} \label{startequation}
  \log (c T) \sum_{d,e \leq N} \frac{a(d)\overline{a(e)}}{[d,e]} 
- \sum_{d,e \leq N} \frac{a(d)\overline{a(e)}}{[d,e]} \cdot \log 
\frac{[d,e]}{(d,e)} 
\end{equation}
with $c > 0$ constant (for example $c = 4 e^{2 \gamma - 1} / 2\pi$). Writing
$
(d,e) = \sum_{\ell | d , \ell | e} \varphi(\ell)
$
diagonalizes the first quadratic form,
\begin{equation} \label{first}
\sum_{d, e \leq N} \frac{a(d)\overline{a(e)}}{[d,e]} = \sum_{\ell} 
\frac{\varphi(\ell)}{\ell^2} \cdot |y(\ell)|^2 
\end{equation}
where
$$
y(\ell) := \sum_{d \leq N} \frac{a(d \ell)}{d}
$$
By Moebius inversion 
$$
1 = \sum_{\ell \leq N} \frac{y(\ell)\mu(\ell)}{\ell}
$$
Thus, by Cauchy-Schwarz,
$$
1 \leq \bigg ( \sum_{n \leq N} \frac{\mu(n)^2}{\varphi(n)} \bigg)^{1/2}
\cdot \bigg ( \sum_{n \leq N} \frac{\varphi(\ell)}{\ell^2} 
\cdot |y(\ell)|^2 \bigg )^{1/2}
$$
It follows that the minimum of the quadratic form
(\ref{first}) is $1 / G$, where
$$
G := \sum_{n \leq N} \frac{\mu(n)^2}{\varphi(n)}
$$
The minimum is attained when $y(\ell) = z(\ell)$ with
$$
z(\ell) = \frac{\mu(\ell)}{G} \cdot \frac{\ell}{\varphi(\ell)} 
$$
The above discussion is subsumed in the lemma below.
\begin{lemma} \label{Lemma1}
With notation as above, we have,
$$
\sum_{d,e \leq N} \frac{a(d) \overline{a(e)}}{[d,e]} 
= \frac{1}{G} + \sum_{\ell \leq N} \frac{\varphi(\ell)}{\ell^2} \cdot |y(\ell) - z(\ell)|^2
$$
\end{lemma}
\begin{Proof}
Expanding the square and using (\ref{first}),
\begin{align*}
\sum_{\ell \leq N} \frac{\varphi(\ell)}{\ell^2} \cdot |y(\ell) - z(\ell)|^2
&  = \sum_{\ell \leq N} \frac{\varphi(\ell)}{\ell^2} |y(\ell)|^2 
- 2 \Re \sum_{\ell \leq N} \frac{\varphi(\ell)}{\ell^2} y(\ell)z(\ell)
+ \sum_{\ell \leq N} \frac{\varphi(\ell)}{\ell^2} z(\ell)^2 \\
& = \sum_{d,e \leq N} \frac{a(d)\overline{a(e)}} {[d,e]} 
- \frac{2}{G} \Re \sum_{\ell \leq N} \frac{\mu(\ell)}{\ell} y(\ell)
+ \frac{1}{G^2} \sum_{\ell \leq N} \frac{\mu(\ell)^2}{\varphi(\ell)} \\
& = \sum_{d,e \leq N} \frac{a(d)\overline{a(e)}} {[d,e]}
 - \frac{2}{G} + \frac{1}{G}
\end{align*}
and the claim follows.
\end{Proof}

It remains to understand the second quadratic form appearing in 
equation (\ref{first}). This is more difficult and is accomplished in the
lemma below. 
\begin{lemma} \label{Lemma2}
Let $N = T^{\theta}$. Then, for $T$ large,
$$
- \sum_{d, e \leq N} \frac{a(d)\overline{a(e)}}
{[d,e]} \cdot \log \bigg ( \frac{[d,e]}{(d,e)} \bigg )
  \geq 1 - (\theta + \varepsilon) \log T \sum_{\ell \leq N}
\frac{\varphi(\ell)}{\ell^2} \cdot |y(\ell) - z(\ell)|^2
$$
\end{lemma}
We will prove Lemma \ref{Lemma2} later on.
Assuming the result of Lemma \ref{Lemma2} Proposition B follows
immediately.
\\
\begin{ProofB}
Take $N = T^{\theta}$ with $\theta < 1$. Let $c > 0$ be a constant. 
By Lemma \ref{Lemma1}, and using $G \sim \log N$, we obtain
\begin{equation} \label{firsteqq}
\log (c T) \sum_{d,e \leq N} \frac{a(d) \overline{a(e)}}{[d,e]} = 
 \frac{1}{\theta} + \log (c T)\sum_{\ell \leq N} \frac{\varphi(\ell)}{\ell^2} \cdot |y(\ell) - z(\ell)|^2 + o(1)
\end{equation}
as $T \rightarrow \infty$. 
By Lemma \ref{Lemma2},
\begin{equation} \label{secondeqq}
- \sum_{d, e \leq N} \frac{a(d)\overline{a(e)}}
{[d,e]} \cdot \log \bigg ( \frac{[d,e]}{(d,e)} \bigg ) 
  \geq 1 - (\theta + \varepsilon) \log T \sum_{\ell \leq N}
\frac{\varphi(\ell)}{\ell^2} \cdot |y(\ell) - z(\ell)|^2
\end{equation}
Adding the equations (\ref{firsteqq}) and (\ref{secondeqq}), we obtain
\begin{align*}
\sum_{d,e \leq N} \frac{a(d) \overline{a(e)}} {[d,e]} & \cdot \log
\bigg ( \frac{c T (d,e)}{[d,e]} \bigg ) \\ & \geq 
1 + \frac{1}{\theta} + (1 - \theta - \varepsilon)
\log T \sum_{\ell \leq N} \frac{\varphi(\ell)}{\ell^2} \cdot
|y(\ell) - z(\ell)|^2 + o(1)
\end{align*}
The sum of squares is non-negative, and $1 - \theta - \varepsilon > 0$;
we thus obtain the desired lower bound
$1 + 1 / \theta + o(1)$.
\end{ProofB}
 
\subsection{Proof of Lemma \ref{Lemma2}}

In order to prove Lemma \ref{Lemma2}
we start by expressing the quadratic form (\ref{startequation}) in terms of the sequence
$y(\ell)$. 

\begin{lemma} \label{firststep}
We have, 
\begin{multline} \label{secondsum}
\sum_{d,e \leq N} \frac{a(d) \overline{a(e)}}{[d,e]} \cdot
\log \frac{[d,e]}{(d,e)} = 2\sum_{p^{\alpha} \ell \leq N}
\frac{\log p}{p^{\alpha}} \cdot \frac{\varphi(\ell)}{\ell^2} \cdot
\Re \big ( y(\ell) \overline{y(p^{\alpha} \ell)} \big ) + \\
+ O \bigg ( (\log \log N)^2 \sum_{\ell \leq N}
\frac{\varphi(\ell)}{\ell^2} \cdot |y(\ell) - z(\ell)|^2 + 
\frac{(\log\log N)^2}{\log N} \bigg )
\end{multline}
\end{lemma}
\begin{Proof}
Since $[d\ell, e\ell]/(d\ell, e\ell) = [d,e]/(d,e)$ using
the formula $(d,e) = \sum_{\ell | d, \ell | e} \varphi(\ell)$ we find
\begin{equation} \label{firsttrans}
\sum_{e,d \leq N} \frac{a(d) \overline{a(e)}}{[d,e]} \cdot \log 
\frac{[d,e]}{(d,e)} = \sum_{\ell \leq N} \frac{\varphi(\ell)}{\ell^2}
\sum_{d,e \leq N/\ell} \frac{a(d\ell)\overline{a(e\ell)}}{d e}
\cdot \log \frac{[d,e]}{(d,e)}.
\end{equation}
A prime $p$ divides $[d,e]/(d,e)$ if and only if $|v_p(e) - v_p(d)| \geq 1$
where $v_p(n)$ denotes the $p$-adic valuation of $n$. Therefore,
$$
\log \frac{[d,e]}{(d,e)} = \sum_{\substack{p^\alpha \| e, p^{\beta} \| d \\ |\alpha
- \beta| \geq 1}} \log p
$$
It follows that (\ref{firsttrans}) can be expressed as
\begin{equation}\label{secondeqqq}
\sum_{|\alpha - \beta| \geq 1} 
\sum_{p \ell \leq N} \frac{\log p}{p^{\alpha + \beta}}\cdot \frac{\varphi(\ell)}
{\ell^2} \cdot \big ( y(p^{\alpha} \ell) - \frac{1}{p} y(p^{\alpha + 1}\ell)
\big ) \cdot \big 
( \overline{ y(p^{\beta} \ell) - \frac{1}{p} y(p^{\beta + 1} \ell)} \big )
\end{equation}
We bound the contribution of the terms with $\alpha, \beta \geq 1$: multiplying
out and repeatedly using
the inequality $2a b \leq a^2 + b^2$ we find that,
$$
\sum_{\substack{|\alpha - \beta | \geq 1 \\ \alpha, \beta \geq 1}} 
\frac{1}{p^{\alpha + \beta}}
\cdot \big ( y(p^{\alpha} \ell) - \frac{1}{p} y(p^{\alpha + 1}\ell)
\big ) \cdot \big ( \overline{y(p^{\beta} \ell) - \frac{1}{p} y(p^{\beta + 1} \ell)} \big )
\ll \sum_{\alpha \geq 1} \frac{|y(p^\alpha \ell)|^2}{p^{\alpha}} \cdot \frac{1}{p}
$$
Inserting this back into (\ref{secondeqqq}), bounds the contributions of
the terms with $\alpha , \beta \geq 1$ by
\begin{align} \label{boundagain}
\ll \sum_{\substack{p^\alpha \ell \leq N \\ \alpha \geq 1}} 
\frac{\log p}{p^{\alpha + 1}} \cdot
\frac{\varphi(\ell)}{\ell^2} \cdot |y(p^\alpha \ell)|^2 
&  = \sum_{m \leq N} \bigg ( \sum_{\substack{p^\alpha \ell = m \\ \alpha \geq 1}} 
\frac{\log p }{p^{\alpha + 1}}
\cdot \frac{\varphi(\ell)}{\ell^2} \bigg ) \cdot |y(m)|^2 \\ \nonumber
& \ll (\log \log N)^2 \cdot \sum_{\ell \leq N} \frac{\varphi(\ell)}{\ell^2}
\cdot |y(\ell)|^2 \\ \nonumber
& \ll \frac{(\log\log N)^2}{\log N} + (\log\log N)^2 \sum_{\ell \leq N}
\frac{\varphi(\ell)}{\ell^2} \cdot |y(\ell) - z(\ell)|^2
\end{align}
On the other hand the remaining terms with $\alpha = 0, \beta \geq 1$
and $\beta = 0, \alpha \geq 1$ in (\ref{secondeqqq}) telescope to
$$
2 \sum_{p \ell \leq N} \frac{\log p}{p} \cdot \frac{\varphi(\ell)}{\ell^2}
\cdot \Re \big ( y(\ell)\overline{y(p \ell)} - \frac{1}{p} \cdot |y(p \ell)|^2
\big)
$$
%
To the above sum we add the contribution of the terms with $p^{\alpha} \ell
\leq N$ and $\alpha > 1$. This contribution is estimated by
in the same way as in (\ref{boundagain}) and therefore negligible. This leads
us to a final estimate of
$$
2 \sum_{p^{\alpha } \ell \leq N} \frac{\log p}{p^{\alpha}} \cdot \frac{\varphi(\ell)}{\ell^2} \cdot \Re \big ( y(\ell) \overline{y(p^{\alpha} \ell)} \big ) 
$$
plus the same error as in (\ref{boundagain}).
\end{Proof}

Write
\begin{multline}
y(\ell)\overline{y(p^{\alpha} \ell)} = (y(\ell) - z(\ell)) \cdot 
\overline{(y(p^{\alpha} \ell) - z(p^{\alpha} \ell))} +
\\ + z(\ell) \cdot \overline{(y(p^{\alpha} \ell) - z(p^{\alpha} \ell))} + 
\overline{z(p^{\alpha} \ell)} 
\cdot
(y(\ell) - z(\ell)) + \overline{z(p^{\alpha} \ell)} z(\ell).
\end{multline}
It follows from the above identity and Lemma \ref{firststep}
that 
$$
\sum_{d,e \leq N} \frac{a(d) \overline{a(e)}}{[d,e]} \cdot \log 
\frac{[d,e]}{(d,e)} = S_1 - S_2 + S_3
$$
where
\begin{align*}
S_1 & := 2\sum_{p^{\alpha} \ell \leq N} \frac{\log p}{p^{\alpha}} \cdot \frac{\varphi(\ell)}
{\ell^2} \cdot (y(\ell) - z(\ell))\overline{(y(p^{\alpha} \ell) - z(p^{\alpha}
 \ell))} \\
S_2 & := 2\sum_{p^{\alpha} \ell \leq N} \frac{\log p}{p^{\alpha}} \cdot 
\frac{\varphi(\ell)}
{\ell^2} \cdot \big ( z(\ell) \overline{(y(p^{\alpha} \ell) - z(p^{\alpha} 
\ell))} 
+ \overline{z(p^{\alpha} \ell)} (y(\ell) - z(\ell)) \big ) \\
S_3 & := 2\sum_{p^{\alpha} \ell \leq N} \frac{\log p}{p^{\alpha}} \cdot 
\frac{\varphi(\ell)}
{\ell^2} \cdot z(\ell) \overline{z(p^{\alpha} \ell)} 
\end{align*}
Lemma \ref{Lemma2} follows upon computing $S_1$, $S_2$ and $S_3$ and combining
the resulting estimate. We perform the necessary computations in the three
lemma below. 

\begin{lemma}
We have, 
\begin{equation*}
|S_1| \leq  (\log N + O (\log \log N)) \cdot \sum_{\ell} \frac{\varphi(\ell)}{\ell^2}
\cdot |y(\ell) - z(\ell)|^2
\end{equation*}
\end{lemma}
\begin{Proof}
Applying to $2(y(\ell) - z(\ell))\overline{(y(p^{\alpha} \ell) - z(p^{\alpha}
 \ell))}$ the
inequality $2 |a b| \leq |a|^2  + |b|^2$ we find
\begin{align*} 
|S_1| & \leq \sum_{p^{\alpha} \ell \leq N} \frac{\log p}{p^{\alpha}} 
\cdot \frac{\varphi(\ell)}
{\ell^2} \cdot |y(\ell) - z(\ell)|^2 + \sum_{p^{\alpha} \ell \leq N}
\frac{\log p}{p^{\alpha}} \cdot \frac{\varphi(\ell)}{\ell^2}
\cdot |y(p^{\alpha} \ell) - z(p^{\alpha} \ell)|^2 \\
 & \leq \sum_{\ell \leq N} \frac{\varphi(\ell)}{\ell^2}
\cdot |y(\ell) - z(\ell)|^2 \cdot \log (N / \ell)
+ 
\sum_{m \leq N} \bigg ( \sum_{p^{\alpha} \ell = m} \frac{\log p}{p^{\alpha}}
\cdot \frac{\varphi(\ell)}{\ell^2} \bigg ) \cdot
|y(m) - z(m)|^2
\end{align*}
For $m = p^{\alpha} \ell$ we have
$\varphi(\ell) / \ell = \varphi(m) / m \cdot (1 + O(1 / p))$. Therefore,
$$
\sum_{p^{\alpha} \ell = m} \frac{\log p}{p^{\alpha}} \cdot \frac{\varphi(\ell)}{\ell^2}
= \frac{\varphi(m)}{m^2} \cdot \log m + O \bigg
(\frac{\varphi(m)}{m^2} \cdot \log\log m \bigg ) 
$$
Therefore the sums with $\log m$ cancel out
and we obtain the bound
$$
(\log N + O(\log \log N) ) \cdot 
\sum_{\ell \leq N} \frac{\varphi(\ell)}{\ell^2} \cdot |y(\ell) - z(\ell)|^2 .
$$
as desired.
\end{Proof}
\begin{lemma}
We have
\begin{equation*}
|S_2| \ll \frac{\log \log N}{\sqrt{\log N}}
\cdot \sum_{\ell \leq N} \frac{\varphi(\ell)}{\ell^2} \cdot
|y(\ell) - z(\ell)|^2 
\end{equation*}
\end{lemma}
\begin{Proof}
On the one hand,
\begin{align} \label{firstequation} \nonumber
 2\sum_{p^{\alpha} \ell \leq N} \frac{\log p}{p^{\alpha}} 
\cdot \frac{\varphi(\ell)}{\ell^2}
\cdot z(\ell) \overline{(y(p^{\alpha} \ell) - z(p^{\alpha} \ell))}  = & 
\frac{2}{G} \sum_{p^{\alpha} \ell \leq N}
\frac{\log p}{p^{\alpha}} \frac{\mu(\ell)}{\ell} \cdot
\overline{(y(p^{\alpha} \ell) - z(p^{\alpha} \ell))} \\ \nonumber
= & \frac{2}{G} \sum_{m \leq N} \bigg ( \sum_{p^{\alpha} \ell = m}
\frac{\log p}{p^{\alpha}} \cdot \frac{\mu(\ell)}{\ell} \bigg )
\cdot \overline{(y(m) - z(m))} \\ 
= & - \frac{2}{G} \sum_{m \leq N} \frac{\mu(m) \log m}{m}
\cdot \overline{(y(m) - z(m))}
\end{align}
On the other hand, 
\begin{align*}
2 \sum_{p^{\alpha} \ell \leq N} \frac{\log p}{p^{\alpha}} 
\cdot \frac{\varphi(\ell)}{\ell^2}
\cdot \overline{z(p^{\alpha} \ell)} ( y(\ell) - z(\ell))
= & - \frac{2}{G} \sum_{\substack{p^{\alpha} \ell \leq N \\ (p, \ell) = 1}} 
\frac{\log p}{p^{\alpha}} \cdot \frac{\mu(\ell)}{\ell} \cdot (y( 
\ell) 
- z(\ell)) \\
= & - \frac{2}{G} \sum_{\ell \leq N} \frac{\mu(\ell)}{\ell} \cdot
(y(\ell) - z(\ell)) \sum_{\substack{p^{\alpha} \leq N/ \ell \\ (p , \ell) = 1}}
\frac{\log p}{p^{\alpha}}
\end{align*}
Since
$
\sum_{\substack{p^{\alpha} \leq N/\ell \\ (p , \ell) = 1}} 
\frac{\log p}{p^{\alpha}}
= \log N / \ell + O(\log \log N)
$
and
$$
\frac{1}{G} \sum_{\ell \leq N} \frac{\mu(\ell)}{\ell}
\cdot (y(\ell) - z(\ell)) = 0
$$
the sum simplifies to
\begin{equation} \label{secondequation}
\frac{2}{G} \sum_{\ell \leq N} \frac{\mu(\ell) \log \ell}{\ell} 
\cdot (y(\ell) - z(\ell))
+ O \bigg ( 
\frac{\log\log N}{\sqrt{\log N}} \cdot
\sum_{\ell \leq N} \frac{\varphi(\ell)}{\ell^2} 
\cdot |y(\ell) - z(\ell)|^2 \bigg )
\end{equation}
Adding (\ref{firstequation}) and (\ref{secondequation}) the main terms cancel and
we obtain the bound for $|S_2|$.
\end{Proof}
\begin{lemma}
We have,
$$
S_3 
 = -1 + O \bigg ( \frac{\log\log N}{\log N} \bigg ).
$$
\end{lemma}
\begin{Proof}
Since, for $\ell \leq N$, 
$$
\sum_{\substack{p \leq N / \ell \\ (p, \ell) = 1}} \frac{\log p}{p}
= \log (N/\ell) + O(\log\log N)
$$
We have
\begin{align*}
 2 \sum_{p \ell \leq N} \frac{\log p}{p} \cdot \frac{\varphi(\ell)}{\ell^2}
\cdot z(\ell) z(p \ell) &  = - 
\frac{2}{G^2} \sum_{\substack{p \ell \leq N \\ (p, \ell) = 1}} 
\frac{\log p}{p} \cdot \frac{\mu(\ell)^2}{\varphi(\ell)}
 \\
& = - \frac{2}{G^2} \sum_{\ell \leq N} \frac{\mu(\ell)^2}{\varphi(\ell)}
\cdot \big ( \log (N / \ell) + O(\log\log N) \big ) \\
& = -1 + O \bigg ( \frac{\log \log N}{\log N} \bigg )
\end{align*}
as desired. 
\end{Proof}

\end{document}